\documentclass{amsart}
\usepackage[latin1]{inputenc}

\usepackage{amssymb,amsmath}
\usepackage{amsthm}
\usepackage{paralist}
\usepackage{url}

\theoremstyle{plain}
\newtheorem{theorem}{Theorem}

\newtheorem*{main}{Associativity Lemma}
\newtheorem{lemma}{Lemma}

\begin{document}

\title[classes and constraints 
 stable under compositions with clones]{Function classes and relational constraints 
 stable under compositions with clones}

\author{Miguel Couceiro}
\address{Mathematics Research Unit, University of Luxembourg \\
162A, avenue de la Fa\"{\i}encerie, L-1511 Luxembourg, Luxembourg }
\email{Miguel.Couceiro@uni.lu}

\author{Stephan Foldes}
\address{Institute of Mathematics, Tampere University of Technology\\ 
PL553, 33101 Tampere, Finland}
\email{sf@tut.fi}

\begin{abstract}
The general Galois theory for functions and 
relational constraints over arbitrary sets described in the authors' previous paper is
 refined by imposing algebraic conditions on relations. 
\end{abstract}

\date{15 November, 2008}

\maketitle

\section{Introduction}

In this paper we extend the results obtained in \cite{CF3} by considering more general closure conditions on 
 classes of functions of several variables, and by restricting relational constraints to consist 
of invariant relations. 
In fact, the Theorems 2.1 and 3.2 in \cite{CF3} correspond to Theorems 1 and 3 below, respectively,
in the particular case $\mathcal{C}_1=\mathcal{C}_2=\mathcal{P}$, where $\mathcal{P}$ 
denotes the smallest clone containing only projections.

A preliminary version of the current manuscript appeared as a Rutcor Research Report
RRR-22-2004 available at \url{http://rutcor.rutgers.edu/~rrr/2004.html}.

 \section{Basic notions and preliminary results}

Throughout the paper, let $A$, $B$, $E$ and $G$ be arbitrary nonempty sets. Given a nonnegative integer $m$, the elements of $A^m$ are viewed as unary functions on the von Neumann ordinal $m=\{0,\ldots ,m-1\}$ to $A$.
 
 A \emph{function of several variables on $A$ to $B$} (or simply, \emph{function on $A$ to $B$})
 is a map $f: A^n \rightarrow B$,
 for some positive integer $n$ called the \emph{arity} of $f$.
A \emph{class} of functions on $A$ to $B$ is a subset $\mathcal{F}\subseteq \cup _{n\geq 1}B^{A^n}$.
For a fixed arity $n$, the $n$ different \emph{projection maps} 
${\bf a}=(a_t\mid t\in n)\mapsto a_i$, $i\in n$, are 
also called \emph{variables}.
For $A=B=\{0,1\}$, a function  on $A$ to $B$ is called a \emph{Boolean function}.

If $f$ is an $n$-ary function on $B$ to $E$ and $g_1,\ldots ,g_n$ are all $m$-ary functions 
on $A$ to $B$ then the \emph{composition} $f(g_1,\ldots ,g_n)$ is an $m$-ary function on $A$ to $E$,
 and its value on $(a_1,\ldots ,a_m)\in A^m$ is $f(g_1(a_1,\ldots ,a_m),\ldots ,g_n(a_1,\ldots ,a_m))$.
If $\mathcal{I}\subseteq \cup _{n\geq 1}E^{B^n}$ and $\mathcal{J}\subseteq \cup _{n\geq 1}B^{A^n}$ we define 
the \emph{composition of} $\mathcal{I}$ \emph{with} $\mathcal{J}$,
 denoted $\mathcal{I}\mathcal{J}$, by
 \begin{displaymath}
 \mathcal{I}\mathcal{J}=\{f(g_1,\ldots ,g_n)\mid n,m\geq 1, f\textrm{ $n$-ary in $\mathcal{I}$, }g_1,\ldots ,g_n\textrm{ $m$-ary in $\mathcal{J}$} \}. 
\end{displaymath}
If $\mathcal{I}$ is a singleton, $\mathcal{I}=\{f\}$, then we write $f\mathcal{J}$ for $\{f\}\mathcal{J}$.   
We say that a class $\mathcal{I}$ of functions of several variables is \emph{stable under right}
 (\emph{left}) \emph{composition with $\mathcal{J}$} if, whenever the composition is well defined,
 $\mathcal{I}\mathcal{J}\subseteq \mathcal{I}$ ( $\mathcal{J}\mathcal{I}\subseteq \mathcal{I}$, respectively).
A \emph{clone} on $A$ is a set $\mathcal{C}\subseteq \cup _{n\geq 1}A^{A^n}$
 that contains all projections and satisfies $\mathcal{C} \mathcal{C}\subseteq \mathcal{C}$ 
(or equivalently, $\mathcal{C} \mathcal{C}=\mathcal{C}$). 
Note that if $\mathcal{J}$ is a clone on $A$ (on $B$) and 
$\mathcal{I}\subseteq \cup _{n\geq 1}B^{A^n}$, then
$\mathcal{I}\mathcal{J}\subseteq \mathcal{I}$ if and only if $\mathcal{I}\mathcal{J}=\mathcal{I}$ 
($\mathcal{J}\mathcal{I}\subseteq \mathcal{I}$ if and only if
$\mathcal{J}\mathcal{I}=\mathcal{I}$, respectively). Note that stability under right composition with the clone $\mathcal{P}\subseteq \cup _{n\geq 1}A^{A^n}$ of all projections on $A$ subsumes the operations of identification of variables, permutation of variables and addition of inessential variables.

\begin{main} 
Let $A$, $B$, $E$ and $G$ be arbitrary nonempty sets, and consider function classes
$\mathcal{I}\subseteq \cup _{n\geq 1}G^{E^n}$, $\mathcal{J}\subseteq \cup _{n\geq 1}E^{B^n}$, and
$\mathcal{K}\subseteq \cup _{n\geq 1}B^{A^n}$. The following hold:
\begin{itemize}
\item[(i)] $(\mathcal{I}\mathcal{J})\mathcal{K}\subseteq \mathcal{I}(\mathcal{J}\mathcal{K})$;
\item[(ii)] If $\mathcal{J}$ is stable under right composition with the clone of projections on $B$,
then $(\mathcal{I}\mathcal{J})\mathcal{K}=\mathcal{I}(\mathcal{J}\mathcal{K})$.
\end{itemize} 
\end{main}

\begin{proof}
The inclusion $(i)$ is a direct consequence of the definition of
function class composition. Property $(ii)$ asserts that the converse
inclusion also holds if $\mathcal{J}$ is stable under right composition with
projections. This hypothesis means in particular that all functions
obtained from members of $\mathcal{J}$ by permutation of variables and addition of inessential variables 
are also in $\mathcal{J}$. A typical function in $\mathcal{I}(\mathcal{J}\mathcal{K})$ is of the form
 \begin{displaymath}
f({g}_1(h_{11},\ldots ,h_{1m_1}), \ldots ,{g}_n(h_{n1},\ldots ,h_{nm_n}))
\end{displaymath} 
where $f$ is in $\mathcal{I}$, the $g_i$'s are in $\mathcal{J}$, and the $h_{ij}$'s are in $\mathcal{K}$. By taking
appropriate functions ${g'}_1,\ldots ,{g'}_n$ obtained from $g_1,\ldots ,g_n$ by permutation of variables and
addition of inessential variables, the function above can be
expressed as
 \begin{displaymath}
f({g'}_1(h_{11},\ldots ,h_{1m_1},\ldots ,h_{n1},\ldots ,h_{nm_n}), \ldots ,{g'}_n(h_{11},\ldots ,h_{1m_1},\ldots ,h_{n1},\ldots ,h_{nm_n}))
\end{displaymath} 
which is easily seen to be in $(\mathcal{I}\mathcal{J})\mathcal{K}$ .
\end{proof}

Note that statement $(ii)$ of the Associativity Lemma applies, in particular, if $\mathcal{J}$ is any clone on $E=B$.

Let $\mathcal{F}$ be a set of functions on $A$ to $B$. If $\mathcal{P}$ is the clone 
of all projections on $A$, then $\mathcal{F}\mathcal{P}=\mathcal{F}$ expresses closure under taking minors in \cite{Pi}, or closure under simple variable substitutions in the terminology of \cite{CF3}.
 If $A=B=\{0,1\}$ and $\mathcal{L}_{01}$ is the clone (Post class) 
of constant preserving linear Boolean functions, then $\mathcal{F}\mathcal{L}_{01}=\mathcal{F}$
 is equivalent to closure under substitution of triple sums $x+y+z$ for variables,
 while $\mathcal{L}_{01}\mathcal{F}=\mathcal{F}$ 
is equivalent to closure under taking triple sums of Boolean functions $f+g+h$ (see \cite{CF1}).

An \emph{$m$-ary relation on $A$} is a subset $R$ of $A^m$.
Thus the relation $R$ is a class (set) of unary maps on $m$ to $A$.
 A function $f$ of several variables on $A$ to $A$ is said to \emph{preserve} $R$ if $fR\subseteq R$.

For a class $\mathcal{F}\subseteq \cup _{n\geq 1}A^{A^n}$ of functions on $A$,
 an $m$-ary relation $R$ on $A$ is called an 
$\mathcal{F}$-\emph{invariant} if $\mathcal{F}R\subseteq R$.
 In other words, $R$ is an 
$\mathcal{F}$-invariant if every member of $\mathcal{F}$ preserves $R$.
If two classes of functions $\mathcal{F}$ and $\mathcal{G}$ generate the same clone, then the      
$\mathcal{F}$-invariants are the same as the $\mathcal{G}$-invariants.
(See P$\ddot o$schel \cite{Po1} and \cite{Po2}.)

Observe that we always have $R\subseteq \mathcal{F}R$ if $\mathcal{F}$ contains the projections,
but we can have $R\subseteq \mathcal{F}R$ even if $\mathcal{F}$ contains no projections.
(Take the Boolean triple sum $x_1+x_2+x_3$ as the only member of $\mathcal{F}$.)

For a clone $\mathcal{C}$, the intersection of $m$-ary $\mathcal{C}$-invariants
 is always a $\mathcal{C}$-invariant and 
it is easy to see that,
for an $m$-ary relation $R$, the smallest $\mathcal{C}$-invariant
 containing $R$ in $A^m$ is $\mathcal{C}R$, and it 
 is said to be \emph{generated by} $R$.
(See \cite{Po1} and \cite{Po2}, 
where P$\ddot o$schel denotes $\mathcal{C}R$ by $\Gamma _{\mathcal{C}}(R)$.)

\section{Classes of Functions Definable by Constraints Consisting of Invariant Relations}

Consider arbitrary non empty sets $A$ and $B$. 
An $m$-ary \emph{$A$-to-$B$ constraint} (or simply, $m$-ary \emph{constraint}, 
when the underlying sets are understood from the context)
is a couple $(R,S)$ where $R\subseteq A^m$ and $S\subseteq B^m$. 
The relations $R$ and $S$ are called the \emph{antecedent} and \emph{consequent}, respectively,
 of the relational constraint (Pippenger \cite{Pi}).
Let $\mathcal{C}_1$ and $\mathcal{C}_2$ be clones on $A$ and $B$, respectively.
If $R$ is a $\mathcal{C}_1$-invariant and $S$ is a $\mathcal{C}_2$-invariant, 
we say that $(R,S)$ is a \emph{$(\mathcal{C}_1,\mathcal{C}_2)$-constraint}.  
 A function 
 $f:A^n\longrightarrow B$, $n\geq 1$, is said to \emph{satisfy} an $m$-ary $A$-to-$B$ constraint $(R,S)$ if
 $fR\subseteq S$.

The following result generalizes Lemma 1 in \cite{CF1}:

\begin{lemma} Consider arbitrary nonempty sets $A$ and $B$. Let $f$ be a 
function  on $A$ to $B$ and let $\mathcal{C}$ be a clone on $A$. 
If every function in $f\mathcal{C}$ satisfies an 
$A$-to-$B$ constraint $(R,S)$,
 then $f$ satisfies $(\mathcal{C}R,S)$.
\end{lemma}
\begin{proof}
The assumption means that $(f\mathcal{C})R\subseteq S$. By the Associativity Lemma,
 $(f\mathcal{C})R=f(\mathcal{C}R)$, and thus $f(\mathcal{C}R)\subseteq S$.  
\end{proof}

A class $\mathcal{K}\subseteq \cup _{n\geq 1}B^{A^n}$ of functions  on $A$ to $B$ is said 
to be \emph{locally closed} if for every function $f$  on $A$ to $B$ the following holds: 
if every finite restriction of $f$ (i.e restriction to a finite subset) coincides with a finite restriction of 
some member of $\mathcal{K}$, then $f$ belongs to $\mathcal{K}$.

A class $\mathcal{K}\subseteq \cup _{n\geq 1}B^{A^n}$ of functions  on $A$ to $B$ 
is said to be \emph{definable} by a set $\mathcal{T}$
 of $A$-to-$B$ constraints, if $\mathcal{K}$ is the class of all 
those functions which satisfy every constraint in $\mathcal{T}$.

\begin{theorem}\label{theoremMain0}
Consider arbitrary nonempty sets $A$ and $B$ and let $\mathcal{C}_1$ and $\mathcal{C}_2$
 be clones on $A$ and $B$, respectively.
For any function class $\mathcal{K}\subseteq \cup _{n\geq 1}B^{A^n}$ 
the following conditions are equivalent:
\begin{itemize}
\item[(i)] $\mathcal{K}$ is locally closed and it is stable both under right composition 
 with $\mathcal{C}_1$ and under 
 left composition with $\mathcal{C}_2$;
\item[(ii)] $\mathcal{K}$ is definable by some set of $(\mathcal{C}_1,\mathcal{C}_2)$-constraints.
\end{itemize}
\end{theorem}

\begin{proof} 
To show that $(ii)\Rightarrow (i)$, assume that 
$\mathcal{K}$ is definable by some set $\mathcal{T}$ of $(\mathcal{C}_1,\mathcal{C}_2)$-constraints.
 For every $(R,S)$ in $\mathcal{T}$, we have $\mathcal{K}R\subseteq S$.
Since $R$ is a $\mathcal{C}_1$-invariant, $\mathcal{K}R=\mathcal{K}(\mathcal{C}_1R)$.
By the Associativity Lemma, $\mathcal{K}(\mathcal{C}_1R)=(\mathcal{K}\mathcal{C}_1)R$, and therefore
$(\mathcal{K}\mathcal{C}_1)R=\mathcal{K}R\subseteq S$.
Since this is true for every $(R,S)$ in $\mathcal{T}$ we must have 
$\mathcal{K}\mathcal{C}_1\subseteq \mathcal{K}$.

 For every $(R,S)$ in $\mathcal{T}$, we have $\mathcal{K}R\subseteq S$, and therefore
 $\mathcal{C}_2(\mathcal{K}R)\subseteq \mathcal{C}_2S$. 
 By the Associativity Lemma, 
$(\mathcal{C}_2\mathcal{K})R\subseteq \mathcal{C}_2(\mathcal{K}R)\subseteq \mathcal{C}_2S$,
and $\mathcal{C}_2S=S$ because $S$ is a $\mathcal{C}_2$-invariant. Thus 
$(\mathcal{C}_2\mathcal{K})R\subseteq S$ for every $(R,S)$ in $\mathcal{T}$, and we must have 
$\mathcal{C}_2\mathcal{K}\subseteq \mathcal{K}$.

To see that $\mathcal{K}$ is locally closed, consider $f\not\in \mathcal{K}$, say of arity $n\geq 1$, 
and let $(R,S)$ be an $m$-ary $(\mathcal{C}_1,\mathcal{C}_2)$-constraint that is satisfied by every function $g$
 in $\mathcal{K}$ but not satisfied by $f$. Hence for some ${\bf a}^1,\ldots ,{\bf a}^n$ in $R$, 
 $f({\bf a}^1,\ldots ,{\bf a}^n)\not\in S$ but $g({\bf a}^1,\ldots ,{\bf a}^n)\in S$, for every 
$n$-ary function $g$ in $\mathcal{K}$.
 Thus the restriction of $f$ to the finite set $\{({\bf a}^1(i),\ldots ,{\bf a}^n(i)):i\in m\}$ 
does not coincide with that of any member of $\mathcal{K}$.

To prove $(i)\Rightarrow (ii)$, we show that for every function $g$ not in $\mathcal{K}$,
there is a $(\mathcal{C}_1,\mathcal{C}_2)$-constraint $(R,S)$ which is satisfied
by every member of $\mathcal{K}$ but not satisfied by $g$.
The class $\mathcal{K}$ will then be definable by the set $\mathcal{T}$ of those
 $(\mathcal{C}_1,\mathcal{C}_2)$-constraints that are satisfied
by all members of $\mathcal{K}$.
 
Note that $\mathcal{K}$ is a fortiori stable under right composition 
with the clone containing all projections, that is, $\mathcal{K}$ is closed under 
simple variable substitutions.
We may assume that $\mathcal{K}$ is non empty.
Suppose that $g$ is an $n$-ary function  on $A$ to $B$ 
not in $\mathcal{K}$. Since $\mathcal{K}$ is locally closed,
 there is a finite restriction $g_F$ of $g$ to a finite subset $F\subseteq A^n$
 such that $g_F$ disagrees with every function in $\mathcal{K}$ restricted to $F$.
Suppose that $F$ has size $m$, and let ${\bf a}^1,\ldots ,{\bf a}^n$ be $m$-tuples in $A^m$, such that
 $F=\{({\bf a}^1(i),\ldots ,{\bf a}^n(i)):i\in m\}$. 
Define $R_0$ to be the set containing ${\bf a}^1,\ldots ,{\bf a}^n$, and
 let $S=\{f({\bf a}^1,\ldots ,{\bf a}^n):f\in \mathcal{K}, f$ $n$-ary$\}$.
 Clearly, $(R_0,S)$ is not satisfied by $g$, 
and it is not difficult to see that every member of 
$\mathcal{K}$ satisfies $(R_0,S)$. 
As $\mathcal{K}$ is stable under left composition with $\mathcal{C}_2$,
 it follows that $S$ is a $\mathcal{C}_2$-invariant.
 Let $R$ be the $\mathcal{C}_1$-invariant generated by $R_0$, i.e. $R=\mathcal{C}_1R_0$.
By Lemma 1, the constraint $(R,S)$ constitutes indeed the desired separating 
$(\mathcal{C}_1,\mathcal{C}_2)$-constraint.
\end{proof}

This generalizes the characterizations of closed classes of functions given by Pippenger in \cite{Pi}
 as well as in \cite{CF1} and \cite{CF3}
 by considering arbitrary underlying sets, possible infinite, 
and more general closure conditions. 
We obtain as special cases of Theorem 1 the characterizations given in 
Theorem 2.1 of \cite{CF3} and, in the finite case, in Theorem 3.2 of \cite{Pi}, by considering $\mathcal{C}_1=\mathcal{C}_2=\mathcal{P}$, and
$\mathcal{C}_1=\mathcal{U}$ and $\mathcal{C}_2=\mathcal{P}$, respectively, where $\mathcal{U}$ is
 a clone containing only functions 
having at most one essential variable, and $\mathcal{P}$ is the clone 
of all projections. Taking $A=B=\{0,1\}$ and 
 $\mathcal{C}_1=\mathcal{C}_2=\mathcal{L}_{01}$, we get the characterization of
 classes of Boolean functions definable by sets of affine constraints 
given in \cite{CF1}.   

\section{Sets of Invariant Constraints Characterized by Functions of Several Variables}

In order to discuss sets of constraints determined by functions of several variables, we need to recall the following concepts and constructions introduced in \cite{Pi} and \cite{CF3}.

Given maps $f:A\rightarrow B$ and $g:C\rightarrow D$, their composition $g\circ f$ is defined only if
$B=C$. Removing this restriction, the \emph{concatenation} of $f$ and $g$, denoted simply
 $gf$, is defined as the map with domain $f^{-1}[B\cap C]$ and codomain $D$ given by 
$(gf)(a)=g(f(a))$ for all $a\in f^{-1}[B\cap C]$. Clearly, if $B=C$ then $gf=g\circ f$, thus concatenation
subsumes and extends functional composition.

Let $(g_i)_{i\in I}$ be a  family of maps, $g_i:A_i\rightarrow B_i$ such that $A_i\cap A_j=\emptyset $
whenever $i\not=j$. The (\emph{piecewise}) \emph{sum of the family} $(g_i)_{i\in I}$, denoted 
${\Sigma }_{i\in I}g_i$, is the map from ${\cup }_{i\in I}A_i$ to ${\cup }_{i\in I}B_i$ whose restriction to each $A_i$ agrees with $g_i$.
If $I$ is finite, we may use the infix $+$ notation.

For $B\subseteq A$, ${\iota }_{AB}$ denotes the canonical injection (inclusion map) from $B$ to $A$.
Note that the restriction $f\mid _B$ of any map $f:A\rightarrow C$
 to the subset $B$ is given by $f\mid _B$ is the concatenation $f{\iota }_{AB}$.

Let  $=_A$ be the equality relation on a set $A$. The \emph{binary $A$-to-$B$ equality constraint} is simply $(=_A,=_B)$.
A constraint $(R,S)$ is called the \emph{empty constraint} if both antecedent and consequent are empty. For every $m\geq 1$, the constraints $(A^m,B^m)$ are said to be \emph{trivial}. Note that every function on $A$ to $B$ satisfies each of these constraints. 

A constraint $(R,S)$ is said to be a \emph{relaxation} of a constraint $(R_0,S_0)$ if $R\subseteq R_0$ and $S\supseteq S_0$.
Given a non-empty family of constraints $(R,S_j)_{j\in J}$ of the same arity (and antecedent),  the constraint $(R,\cap _{j\in J} S_j)$ is said to be obtained from $(R,S_j)_{j\in J}$ by \emph{intersecting consequents}. 

Let $m$ and $n$ be positive integers (viewed as ordinals, i.e., $m=\{0,\ldots ,m-1\}$).
Let $h:n\rightarrow m\cup V$ where
$V$ is an arbitrary set of symbols disjoint from the ordinals 
called \emph{existentially quantified indeterminate indices}, or simply \emph{indeterminates},
 and $\sigma :V\rightarrow A$ any map called a \emph{Skolem map}.
Then each $m$-tuple ${\bf a}\in A^m$, being a map ${\bf a}:m\rightarrow A$,
gives rise to an $n$-tuple $({\bf a}+\sigma )h\in A^n$.

Let $H=(h_j)_{j\in J}$ be a non-empty family of maps $h_j:n_j\rightarrow m\cup V$, where each $n_j$ is a positive integer
(recall $n_j=\{0,\ldots ,n_j-1\}$). Then $H$ is called a \emph{minor formation scheme} with \emph{target} $m$,
 \emph{indeterminate set} $V$ and \emph{source family} $(n_j)_{j\in J}$.
Let $(R_j)_{j\in J}$ be a family of relations (of various arities) on the same set $A$, each $R_j$ of arity $n_j$,
and let $R$ be an $m$-ary relation on $A$. We say that $R$ is a  
 \emph{restrictive conjunctive minor} of the family $(R_j)_{j\in J}$ \emph{via $H$},
or simply a \emph{restrictive conjunctive  minor} of the family $(R_j)_{j\in J}$, if 
for every $m$-tuple $\bf a$ in $A^m$, the condition $R({\bf a})$  
implies that there is a Skolem map $\sigma :V\rightarrow A$ such that, for all $j$ in $J$, we have
$R_j[({\bf a}+\sigma )h_j]$.
On the other hand, if for every $m$-tuple $\bf a$ in $A^m$, the condition $R({\bf a})$ holds whenever
there is a Skolem map $\sigma :V\rightarrow A$ such that, for all $j$ in $J$, we have 
$R_j[({\bf a}+\sigma )h_j]$, then
we say that $R$ is an \emph{extensive conjunctive minor} of the family $(R_j)_{j\in J}$ \emph{via $H$},
or simply an \emph{extensive conjunctive minor} of the family $(R_j)_{j\in J}$.
If $R$ is both a restrictive conjunctive minor and 
an extensive conjunctive minor of the family $(R_j)_{j\in J}$ via $H$, 
then $R$ is said to be a \emph{tight conjunctive minor} of the family $(R_j)_{j\in J}$ \emph{via $H$},
or \emph{tight conjunctive minor} of the family. Note that given a scheme $H$ and a family $(R_j)_{j\in J}$,
 there is a unique tight conjunctive minor of the family $(R_j)_{j\in J}$ via $H$.

If $(R_j,S_j)_{j\in J}$ is a family of $A$-to-$B$ constraints (of various arities) and $(R,S)$ is an 
$A$-to-$B$ constraint such that for a scheme $H$
\begin{itemize}
\item[(i)] $R$ is a restrictive conjunctive minor of $(R_j)_{j\in J}$ via $H$, 
\item[(ii)] $S$ is an extensive conjunctive minor of $(S_j)_{j\in J}$ via $H$, 
\end{itemize}
then $(R,S)$ is said to be a \emph{conjunctive minor} of the family $(R_j,S_j)_{j\in J}$ \emph{via $H$},
or simply a \emph{conjunctive minor} of the family of constraints.

If both $R$ and $S$ are tight conjunctive minors of the respective families via $H$, the constraint   
$(R,S)$ is said to be a \emph{tight conjunctive minor} of the family $(R_j,S_j)_{j\in J}$ \emph{via $H$}, 
or simply a \emph{tight conjunctive minor} of the family of constraints.
Note that given a scheme $H$ and a family $(R_j,S_j)_{j\in J}$, 
 there is a unique tight conjunctive minor of the family via the scheme $H$. 

We say that a class $\mathcal{T}$ of relational constraints is 
\emph{closed under formation of conjunctive minors}
if whenever every member of the nonempty family $(R_j,S_j)_{j\in J}$ of constraints
 is in $\mathcal{T}$, all conjunctive minors of the family $(R_j,S_j)_{j\in J}$ are also in $\mathcal{T}$.  

The following lemma was first obtained in \cite{CF3} and it shows that closure under 
 formation of conjunctive minors is a necessary condition to describe those sets of constraints determined by functions of several variables.

\begin{lemma}\label{lemma1}
Let $(R,S)$ be a conjunctive minor of a non-empty family $(R_j,S_j)_{j\in J}$ of $A$-to-$B$ constraints.
If $f:A^n\rightarrow B$ satisfies every $(R_j,S_j)$ then $f$ satisfies $(R,S)$.
\end{lemma}

A set $\mathcal{T}$ of relational constraints is said to be  \emph{locally closed} if 
for every $A$-to-$B$ constraint $(R,S)$ the following holds: 
if every relaxation of $(R,S)$ with finite antecedent coincides with 
some member of $\mathcal{T}$, then $(R,S)$ belongs to $\mathcal{T}$.
The following result was shown in \cite{CF3} (see Theorem 3.2) and it provides necessary and sufficient conditions for a set of constraints to be determined by functions of several variables.  

\begin{theorem}\label{theorem2}
Consider arbitrary non-empty sets $A$ and $B$.
Let $\mathcal{T}$ be a set of $A$-to-$B$ relational constraints. Then the following are equivalent:
\begin{itemize}
\item[(i)]$\mathcal{T}$ is locally closed and contains the binary equality constraint,
 the empty constraint, and it is closed under formation of conjunctive minors;
\item[(ii)]  There is a set of functions on $A$ to $B$ which satisfy exactly those constraints in $\mathcal{T}$.
\end{itemize}
\end{theorem}

Let $\mathcal{C}_1$ and $\mathcal{C}_2$ be clones on arbitrary nonempty sets $A$ and $B$, respectively. 
Among all $A$-to-$B$ constraints, observe that the empty constraint and the equality constraint are 
$(\mathcal{C}_1,\mathcal{C}_2)$-constraints.

The following Lemma is essentially a restatement, in a variant form, of
the closure condition given by Szab$\acute o$ in \cite{Sz} on the set of relations
preserved by a clone of functions. We indicate a proof via Lemma~\ref{lemma1} above.

\begin{lemma}\emph{(Szab$\acute o$)}  Let $\mathcal{C}$ be a clone on an arbitrary nonempty set $A$. 
If $R$ is a tight conjunctive minor of a nonempty family $(R_j)_{j\in J}$ of $\mathcal{C}$-invariants,
then $R$ is a $\mathcal{C}$-invariant.
\end{lemma}

\begin{proof} Let $R$ be a tight conjunctive minor of a nonempty family $(R_j)_{j\in J}$ of 
$\mathcal{C}$-invariants. We have to prove that every function in
$\mathcal{C}$ preserves $R$ or, equivalently, that every function in
$\mathcal{C}$ satisfies the $A$-to-$A$ constraint $(R,R)$.
 Since $(R_j)_{j\in J}$ is a nonempty family of 
$\mathcal{C}$-invariants, every function in
$\mathcal{C}$ preserves every member of the family $(R_j)_{j\in J}$, that is, every function in
$\mathcal{C}$ satisfies every member of the family $(R_j,R_j)_{j\in J}$ of $A$-to-$A$ constraints. 
From Lemma~\ref{lemma1} above, it follows that every member of
$\mathcal{C}$ satisfies $(R,R)$, that is, $R$ is a $\mathcal{C}$-invariant.
\end{proof}

Thus every tight conjunctive minor $(R,S)$ of 
a nonempty family 
$(R_j,S_j)_{j\in J}$ of
$(\mathcal{C}_1,\mathcal{C}_2)$-constraints 
 is a $(\mathcal{C}_1,\mathcal{C}_2)$-constraint.
However, not all relaxations of $(\mathcal{C}_1,\mathcal{C}_2)$-constraints 
 are $(\mathcal{C}_1,\mathcal{C}_2)$-constraints and so not all conjunctive minors of 
a nonempty family 
$(R_j,S_j)_{j\in J}$ of
$(\mathcal{C}_1,\mathcal{C}_2)$-constraints 
 are $(\mathcal{C}_1,\mathcal{C}_2)$-constraints.
A relaxation $(R,S)$ of an $A$-to-$B$ constraint $(R_0,S_0)$ 
is called a \emph{$(\mathcal{C}_1,\mathcal{C}_2)$-relaxation} of $(R_0,S_0)$ if $(R,S)$
 is a $(\mathcal{C}_1,\mathcal{C}_2)$-constraint.
Similarly, a conjunctive minor $(R,S)$ of a nonempty family 
$(R_j,S_j)_{j\in J}$ of $A$-to-$B$ constraints 
is called a \emph{$(\mathcal{C}_1,\mathcal{C}_2)$-conjunctive minor} of the family 
$(R_j,S_j)_{j\in J}$, if $(R,S)$
 is a $(\mathcal{C}_1,\mathcal{C}_2)$-constraint.

A set $\mathcal{T}$ of $(\mathcal{C}_1,\mathcal{C}_2)$-constraints is said to be
 \emph{closed under formation of $(\mathcal{C}_1,\mathcal{C}_2)$-conjunctive minors}
 if whenever every member of the nonempty family $(R_j,S_j)_{j\in J}$ of constraints
 is in $\mathcal{T}$, all $(\mathcal{C}_1,\mathcal{C}_2)$-conjunctive minors of the family $(R_j,S_j)_{j\in J}$ are also in $\mathcal{T}$.
The following result extends Lemma 1 in \cite{CF2}.

\begin{lemma}\label{lemma3} Let $\mathcal{C}_1$ and $\mathcal{C}_2$ be clones on arbitrary nonempty sets $A$ and $B$,
 respectively. Let $\mathcal{T}_0$ be a set 
of $(\mathcal{C}_1,\mathcal{C}_2)$-constraints, closed under $(\mathcal{C}_1,\mathcal{C}_2)$-relaxations.
 Define $\mathcal{T}$ 
to be the set of all relaxations of the various constraints in $\mathcal{T}_0$.
 Then $\mathcal{T}_0$ is the set of $(\mathcal{C}_1,\mathcal{C}_2)$-constraints which are in
 $\mathcal{T}$, and the following are equivalent:
\begin{itemize}
\item[(a)]$\mathcal{T}_0$ is
closed under formation of $(\mathcal{C}_1,\mathcal{C}_2)$-conjunctive minors;
\item[(b)] $\mathcal{T}$ is closed under taking conjunctive minors.
\end{itemize}
\end{lemma}

\begin{proof} Clearly, the first claim holds, and it is easy to see that $(b)\Rightarrow (a)$. To prove implication $(a)\Rightarrow (b)$, assume $(a)$. 
Let $(R,S)$ be a conjunctive minor of a nonempty family 
$(R_j,S_j)_{j\in J}$ of $A$-to-$B$ constraints in $\mathcal{T}$ via a scheme $H=(h_j)_{j\in J}$, 
$h_j:n_j\rightarrow m\cup V$.
We have to prove that $(R,S)\in \mathcal{T}$.
 
Since for every $j$ in $J$ $(R_j,S_j)\in \mathcal{T}$, there is a nonempty family 
$({R}_j^0,{S}_j^0)_{j\in J}$ of $(\mathcal{C}_1,\mathcal{C}_2)$-constraints 
in $\mathcal{T}_0$ such that, for each $j$ in $J$, $(R_j,S_j)$ is a relaxation of $({R}_j^0,{S}_j^0)$.
 So let $({R}_0,{S}_0)$ be the tight conjunctive minor of the family $({R}_j^0,{S}_j^0))_{j\in J}$ 
 via the scheme $H$. From Lemma 2, it follows that $R_0$ is a $\mathcal{C}_1$-invariant
and $S_0$ a $\mathcal{C}_2$-invariant, and
since $\mathcal{T}_0$ is closed under formation of
$(\mathcal{C}_1,\mathcal{C}_2)$-conjunctive
 minors, we have $({R}_0,{S}_0)\in \mathcal{T}_0$.

Let us prove that $(R,S)$ is a relaxation of $({R}_0,{S}_0)$ and, thus, that $(R,S)\in \mathcal{T}$.
Since $R$ is a restrictive conjunctive minor of the family 
$(R_j)_{j\in J}$ via the scheme $H=(h_j)_{j\in J}$, we have that for every
 $m$-tuple ${\bf a}$ in $R$ there is a Skolem map 
$\sigma :V\rightarrow A$ such that, for all $j$ in $J$,
the $n_j$-tuple $({\bf a}+\sigma )h_j$ is in $R_j$. 
 Since $R_j\subseteq {R}_j^0$ for every $j$ in $J$, it follows that $({\bf a}+\sigma )h_j$ is in ${R}_j^0$
for every $j$ in $J$.
Thus ${\bf a}$ is in $R_0$ and we conclude $R\subseteq R_0$. 

By analogous reasoning one can easily verify that ${\bf b}$ is in $S$ whenever ${\bf b}$ is in $S_0$, i.e
that $S\supseteq S_0$. Thus $(R,S)$ is a relaxation of $(R_0,S_0)$
 and so $(R,S)\in \mathcal{T}$, and the proof of \emph{(a)} is complete.
\end{proof}

Let $\mathcal{T}_0$ be a set of 
$(\mathcal{C}_1,\mathcal{C}_2)$-constraints. We say that $\mathcal{T}_0$ is 
$(\mathcal{C}_1,\mathcal{C}_2)$-\emph{locally closed} if 
 the set $\mathcal{T}$ of all relaxations of the various constraints in $\mathcal{T}_0$ is locally closed.

We can now extend Theorem~\ref{theorem2} above to sets of $(\mathcal{C}_1,\mathcal{C}_2)$-constraints.

\begin{theorem}\label{theoremMain} Let $\mathcal{C}_1$ and $\mathcal{C}_2$ be clones on arbitrary nonempty sets $A$ and $B$,
 respectively, and let $\mathcal{T}_0$ be a set of $(\mathcal{C}_1,\mathcal{C}_2)$-constraints.
 Then the following are equivalent:
\begin{itemize}
\item[(i)]$\mathcal{T}_0$ is $(\mathcal{C}_1,\mathcal{C}_2)$-locally closed,
 contains the binary equality constraint,
 the empty constraint, and it is closed under formation of
$(\mathcal{C}_1,\mathcal{C}_2)$-conjunctive minors;
\item[(ii)] There is a set of functions  on $A$ to $B$ which satisfy exactly those 
$(\mathcal{C}_1,\mathcal{C}_2)$-constraints that are in $\mathcal{T}_0$.
\end{itemize}
\end{theorem}

\begin{proof} To prove implication \emph{(ii) $\Rightarrow $(i)}, assume $(ii)$.
Let $\mathcal{K}$ be the set of all functions 
satisfying every constraint in $\mathcal{T}_0$. Note that $\mathcal{T}_0$ is
closed under $(\mathcal{C}_1,\mathcal{C}_2)$-relaxations. By Theorem 1, we have 
$\mathcal{C}_2\mathcal{K}=\mathcal{K}$, and $\mathcal{K}\mathcal{C}_1=\mathcal{K}$.
 We may assume that $\mathcal{K}\not=\emptyset $.
 Let $\mathcal{T}$ be the set 
 of all those constraints (not necessarily $(\mathcal{C}_1,\mathcal{C}_2)$-constraints) 
satisfied by every function in $\mathcal{K}$. Observe that $\mathcal{T}_0$ is the set of all
$(\mathcal{C}_1,\mathcal{C}_2)$-constraints which are in $\mathcal{T}$.
 We show that $\mathcal{T}$ is the
set of all relaxations in $\mathcal{T}_0$. 

Let $(R,S)$ be a constraint in $\mathcal{T}$. From the definition of $\mathcal{T}$, it follows that
$\mathcal{K}R\subseteq S$. Note that 
$\mathcal{K}$ is stable under right composition with the clone of projections on $A$, because 
$\mathcal{K}\mathcal{C}_1=\mathcal{K}$. Thus by the Associativity Lemma it follows that
$\mathcal{C}_2(\mathcal{K}R)=(\mathcal{C}_2\mathcal{K})R$.
 Since $\mathcal{C}_2\mathcal{K}=\mathcal{K}$, we have that 
$\mathcal{C}_2(\mathcal{K}R)=\mathcal{K}R$, i.e. $\mathcal{K}R$
is a $\mathcal{C}_2$-invariant. Also, again because $\mathcal{K}\mathcal{C}_1=\mathcal{K}$, by Lemma 1 we
conclude that every function in $\mathcal{K}$ satisfies $(\mathcal{C}_1R,\mathcal{K}R)$.   
Clearly, $(\mathcal{C}_1R,\mathcal{K}R)$ is a $(\mathcal{C}_1,\mathcal{C}_2)$-constraint,
 therefore it belongs to $\mathcal{T}_0$.
 Thus every constraint $(R,S)$ in $\mathcal{T}$ is a relaxation of a member of
$\mathcal{T}_0$, namely, a relaxation of $(\mathcal{C}_1R,\mathcal{K}R)$.    

 By Theorem~\ref{theorem2} above, we have that $\mathcal{T}$
 is locally closed and contains the binary equality constraint,
 the empty constraint, and it is closed under formation of conjunctive minors.
Since the binary equality constraint and the empty constraint are 
$(\mathcal{C}_1,\mathcal{C}_2)$-constraints, it follows from Lemma~\ref{lemma3} that $(i)$ holds.

To prove implication \emph{(i) $\Rightarrow $(ii)}, it is enough to show that 
for every $(\mathcal{C}_1,\mathcal{C}_2)$-constraint $(R,S)$ not in $\mathcal{T}_0$,
 there is a function $g$ which satisfies every constraint in $\mathcal{T}_0$, but does not satisfy $(R,S)$.

Let $\mathcal{T}$ be the set of relaxations of the various 
$(\mathcal{C}_1,\mathcal{C}_2)$-constraints in $\mathcal{T}_0$.
Observe that $(R,S)\not\in \mathcal{T}$, otherwise $(R,S)$ would be a 
$(\mathcal{C}_1,\mathcal{C}_2)$-relaxation of
 some $(\mathcal{C}_1,\mathcal{C}_2)$-constraint in $\mathcal{T}_0$,
 contradicting the fact implied by $(i)$ that $\mathcal{T}_0$ is closed under taking 
$(\mathcal{C}_1,\mathcal{C}_2)$-relaxations. Clearly, $\mathcal{T}$ is locally closed,
 contains the binary equality constraint, and the empty constraint.
 From Lemma~\ref{lemma3}, it follows that $\mathcal{T}$ is closed under taking conjunctive minors. 
By Theorem~\ref{theorem2},
 there is a function $g$ which does not satisfy $(R,S)$ 
but satisfies every constraint in $\mathcal{T}$ and so, in particular,
 $g$ satisfies every constraint in $\mathcal{T}_0$.
 Thus we have $(i)\Rightarrow (ii)$.
\end{proof}

Theorem~\ref{theoremMain} generalizes the characterizations of closed classes of constraints given in Pippenger \cite{Pi}
 and also in \cite{CF2} as well as \cite{CF3} by considering both arbitrary, possibly infinite, underlying sets, and
more general closure conditions on 
classes of relational constraints. 

Theorems~\ref{theoremMain0} and \ref{theoremMain} may also be viewed as analogues, with constraints instead of relations,
 of the characterization given by P$\ddot o$schel, as part of Theorem 3.2 in \cite{Po3}, 
of the closed sets in a class of Galois connections between operations and relations of a prescribed type on a set $A$.

\end{document}